\documentclass[a4paper,11pt]{amsart}
\usepackage{graphicx}
\usepackage{mathptmx}
\usepackage{amsmath}
\usepackage{amssymb}
\usepackage{enumitem}
\usepackage{xcolor}

\newmuskip\pFqmuskip

\newcommand*\pFq[6][8]{%
  \begingroup 
  \pFqmuskip=#1mu\relax
  \mathcode`=\string"8000
  \begingroup\lccode`\~=`\,
  \lowercase{\endgroup\let~}\pFqcomma
  F^{#2}_{#3}{\left(\genfrac..{0pt}{}{#4}{#5}\bigg|#6\right)}%
  \endgroup
}
\newcommand{\pFqcomma}{\mskip\pFqmuskip}

\newtheorem{theorem}{Theorem}
\newtheorem{lemma}[theorem]{Lemma}
\newtheorem{corollary}[theorem]{Corollary}
\newtheorem{proposition}[theorem]{Proposition}

\begin{document}

\title[$r$-extended Lah-Bell numbers and polynomials associated with $r$-Lah numbers]{$r$-extended Lah-Bell numbers and polynomials associated with $r$-Lah numbers}

\author{Dae San Kim}
\address{Department of Mathematics, Sogang University, Seoul 121-742, Republic of Korea}
\email{dskim@sogang.ac.kr}

\author{Taekyun  Kim}
\address{Department of Mathematics, Kwangwoon University, Seoul 139-701, Republic of Korea}
\email{tkkim@kw.ac.kr}


\begin{abstract}
The aim of this paper is to study some basic properties of the $r$-extended Lah-Bell numbers and polynomials associated with $r$-Lah numbers and to show the connection between the $r$-extended Lah-Bell polynomials and the rising factorial moments of the Poisson random variable shifted by an even nonnegative integer $2r$. Here the nth $r$-extended Lah-Bell number counts the number of ways a set of $n+r$ elements can be partitioned into non-empty linearly ordered subsets such that $r$ distinguished elements have to be in distinct linearly ordered subsets and the $r$-extended Lah-Bell polynomials are natural extensions of the $r$-extended Lah-Bell numbers.
\end{abstract}

\maketitle
 
\section{Introduction and preliminaries}

Throughout this paper, {\it{a non-empty linearly ordered subset will be simply called an ordered block}}. 
Let $n,k,r$ be nonnegative integers, with $n\ge k$.
Then, the $r$-Lah number $L_{r}(n,k)$ counts the number of partitions of a set with $n+r$ elements into $k+r$ ordered blocks such that $r$ distinguished elements have to be in distinct ordered blocks (see [12]). It is natural and meaningful to define {\it{the $r$-extended Lah-Bell number $B_{n,r}^{L}$ as the number of ways a set of $n+r$ elements can be partitioned into ordered blocks such that $r$ distinguished elements have to be in distinct ordered blocks}}. 
From the definitions of $r$-Lah numbers and the $r$-extended Lah-Bell numbers, we note that 
\begin{equation}\label{-1}
B_{n,r}^{L}=\sum_{k=0}^{n}L_{r}(n,k),\quad(n\ge 0). 
\end{equation} 
The $r$-extended Lah-Bell polynomials $B_{n,r}^{L}(x)$ are also defined as natural extensions of the $r$-extended Lah-Bell numbers. For $r=0$, we note that $L(n,k)=L_{0}(n,k),\,B_{n}^{L}=B_{n,0}^{L}=\sum_{k=0}^{n}L(n,k),\,B_{n}^{L}(x)=B_{n,0}^{L}(x)$ are respectively the Lah numbers, the Lah-Bell numbers and Lah-Bell polynomials (see [4]). \par
The aim of this paper is to study some basic properties of the $r$-extended Lah-Bell numbers $B_{n,r}^{L}$ and polynomials $B_{n,r}^{L}(x)$ associated with $r$-Lah numbers $L_{r}(n,k)$ and to show the connection between the $r$-extended Lah-Bell polynomials $B_{n,r}^{L}(x)$ and the rising factorial moments of $X+2r$, where $X$ is a Poisson random variable and $r$ is a nonnegative integer.  \par
The outline of our main results is as follows. We obtain a recurrence relation for $r$-Lah numbers and the generating function of the $r$-extended Lah-Bell numbers. Also, we derive Dobinski-like formulas for the $r$-extended Lah-Bell numbers and polynomials. We express the $r$-extended Lah-Bell numbers in terms of the 
$r$-extended Bell numbers and vice versa. Also, we get an expression of $r$-Lah numbers in terms of the Stirling polynomials of the second kind and vice versa. In addition, we deduce a relation among $r$-Lah numbers. Lastly, as an application, we show that the rising factorial moments of $X+2r$ is equal to the $r$-extended Lah-Bell numbers evaluated at $\alpha$, where $X$ is the Poisson random variable with parameter $\alpha (>0)$, and $r$ is a nonnegative integer.\par
The novelty of this paper is that it reveals the connection between the rising factorial moments
of the Poisson random variable shifted by an even nonnegative integer $2r$ and the $r$-extended Lah-Bell polynomials. This follows from the fact that the probability-generating function of the Poisson random variable shifted by an even nonnegative integer $2r$ is equal to the generating function of $r$-exteneded Lah-Bell polynomials. For the rest of this section, we recall the necessary facts that will be needed throughout this paper.

\vspace{0.1in}

For any integers $n,k\ge 0$, the Lah numbers are given by 
\begin{equation}\label{0}
L(n,k)=\binom{n-1}{k-1}\frac{n!}{k!},\quad(\mathrm{see}\ [3,4,12]). 
\end{equation}
Recently, the $n$th Lah-Bell number $B_{n}^{L}$, $(n\ge 0)$, is defined as the number of ways a set of $n$ elements can be partitioned into ordered blocks. Thus, we have 
\begin{equation}
	B_{n}^{L}=\sum_{k=0}^{n}L(n,k),\quad(n\ge 0),\quad(\mathrm{see}\ [4]). \label{1}
\end{equation}
The Lah-Bell polynomials are given by 
\begin{equation}
	B_{n}^{L}(x)=\sum_{k=0}^{n}L(n,k)x^{k},\quad(\mathrm{see}\ [4]). \label{2}
\end{equation} 
The generating function of Lah numbers is given by 
\begin{equation}
	\frac{1}{k!}\bigg(\frac{1}{1-t}-1\bigg)^{k}=\sum_{n=k}^{\infty}L(n,k)\frac{t^{n}}{n!},\quad(k\ge 0),\quad(\mathrm{see}\ [4]), \label{3}
\end{equation}
and that of Lah-Bell polynomials is given by 
\begin{equation}
	e^{x\big(\frac{1}{1-t}-1\big)}=\sum_{n=0}^{\infty}B_{n}^{L}(x)\frac{t^{n}}{n!},\quad(\mathrm{see}\ [4]). \label{4}
\end{equation}
The falling factorial sequence is defined by 
\begin{displaymath}
	(x)_{0}=1,\quad (x)_{n}=x(x-1)\cdots(x-n+1),\quad(n\ge 1).
\end{displaymath}
The Stirling numbers of the first kind are defined as 
\begin{equation}
	(x)_{n}=\sum_{k=0}^{n}S_{1}(n,k)x^{k},\quad(n\ge 0),\quad(\mathrm{see}\ [1-12]). \label{5}
\end{equation}
As an inversion formula of \eqref{5}, the Stirling numbers of the second kind are given by 
\begin{equation}
	x^{n}=\sum_{k=0}^{n}S_{2}(n,k)(x)_k. \label{6}
\end{equation}
From \eqref{5} and \eqref{6}, we note that the generating function of Stirling numbers of the first kind and that of the second kind are respectively given by
\begin{equation}
	\frac{1}{k!}\big(\log(1+t)\big)^{k}=\sum_{n=k}^{\infty}S_{1}(n,k)\frac{t^{n}}{n!}, \label{7}
\end{equation}
and 
\begin{displaymath}
	\frac{1}{k!}\big(e^{t}-1\big)^{k}=\sum_{n=k}^{\infty}S_{2}(n,k)\frac{t^{n}}{n!},\quad(k\ge 0),\quad(\mathrm{see}\ [1-12]). 
\end{displaymath}
It is known that the Stirling polynomials of the second kind are defined by 
\begin{equation}
	\frac{1}{k!}e^{xt}\big(e^{t}-1\big)^{k}=\sum_{n=k}^{\infty}S_{2}(n,k|x)\frac{t^{n}}{n!},\quad(k\ge 0),\quad(\mathrm{see}\ [9]). \label{8}
\end{equation}
The Bell polynomials are given by
\begin{equation}
	e^{x(e^{t}-1)}=\sum_{n=0}^{\infty}B_{n}(x)\frac{t^{n}}{n!},\quad(\mathrm{see}\ [9]). \label{9}
\end{equation}
For $x=1$, $B_{n}=B_{n}(1)$, $(n\ge 0)$, are called the Bell numbers. \par 

From \eqref{8} and \eqref{9}, we have 
\begin{equation}
	B_{n}(x)=\sum_{k=0}^{n}S_{2}(n,k)x^{k},\quad (n\ge 0),\quad(\mathrm{see}\ [9]). \label{10}
\end{equation}
A random variable $X$, taking on one of the values $0,1,2,\dots$, is said to be the Poisson random variable with parameter $\alpha(>0$), if the probability mass function of $X$ is given by 
\begin{displaymath}
	p(i)=P\{X=i\}=e^{-\alpha}\frac{\alpha^{i}}{i!},\quad i=0,1,2,\dots. 
\end{displaymath}
Let $f(x)$ be a real valued function, and let $X$ be a Poisson random variable. Then we have 
\begin{equation}
	E[f(X)]=\sum_{i=0}^{\infty}f(i)p(i), \quad(\mathrm{see}\ [13]). \label{11}
\end{equation}
For $n,k\ge 0$, with $n\ge k$ and $r\ge 0$, the $r$-Lah number $L_{r}(n,k)$ counts the number of partitions of a set with $n+r$ elements into $k+r$ ordered blocks such that $r$ distinguished elements have to be in distinct ordered blocks. \par 
It is known that the $r$-Lah numbers are given by 
\begin{equation}
	L_{r}(n,k)=\frac{n!}{k!}\binom{n+2r-1}{k+2r-1},\quad(\mathrm{see}\ [12, \mathrm{Theorem}\,\, 3.7]), \label{12}
\end{equation}
From \eqref{12}, we note that the generating function of $r$-Lah numbers is given by 
\begin{equation}
	\sum_{n=k}^{\infty}L_{r}(n,k)\frac{t^{n}}{n!}=\frac{1}{k!}\bigg(\frac{1}{1-t}-1\bigg)^{k}\bigg(\frac{1}{1-t}\bigg)^{2r},\quad(\mathrm{see}\ [12]). \label{13}
\end{equation}

\section{$r$-Lah numbers and $r$-extended Lah-Bell numbers and polynomials} 
Let $n,k,r$ be nonnegative integers with $r\ge 0$. First, we recall that the $r$-extended Lah-Bell number $B_{n,r}^{L}$ is defined as the number of ways a set of $n+r$ elements can be partitioned into ordered blocks such that $r$ distinguished elements have to be in distinct ordered blocks and from \eqref{-1} that we have

\begin{equation*}
B_{n,r}^{L}=\sum_{k=0}^{n}L_{r}(n,k),\quad(n\ge 0). \label{14} 
\end{equation*} \par
On the one hand, we have
\begin{align}
\bigg(\frac{1}{1-t}\bigg)^{x}\bigg(\frac{1}{1-t}\bigg)^{2r}\ &=\ \bigg(1+\frac{t}{1-t}\bigg)^{x}\bigg(\frac{1}{1-t}\bigg)^{2r} \label{16} \\	
&=\ \sum_{k=0}^{\infty}(x)_{k}\frac{1}{k!}\bigg(\frac{t}{1-t}\bigg)^{k}\bigg(\frac{1}{1-t}\bigg)^{2r} \nonumber \\
&=\ \sum_{n=0}^{\infty}\bigg(\sum_{k=0}^{n}L_{r}(n,k)(x)_{k}\bigg)\frac{t^{n}}{n!}. \nonumber
\end{align}
On the other hand, we have
\begin{equation}
	\bigg(\frac{1}{1-t}\bigg)^{x}\bigg(\frac{1}{1-t}\bigg)^{2r}=\bigg(\frac{1}{1-t}\bigg)^{x+2r}=\sum_{n=0}^{\infty}\langle x+2r\rangle_{n}\frac{t^{n}}{n!},\label{17}
\end{equation}
where $\langle x\rangle_{0}=1$, $\langle x\rangle_{n}=x(x+1)\cdots(x+n-1)$, $(n\ge 1)$. \par 
Therefore, by \eqref{16} and \eqref{17}, we obtain the following lemma. 
\begin{lemma}
	For $n,r\ge 0$, we have 
	\begin{displaymath}
		\langle x+2r\rangle_{n}=\sum_{k=0}^{n}L_{r}(n,k)(x)_{k}.
	\end{displaymath}
\end{lemma}
From Lemma 1, we obtain
\begin{align}
	\langle x+2r\rangle_{n+1}\ &=\ \langle x+2r\rangle_{n}(x+2r+n) \label{18} \\
	&=\ \sum_{k=0}^{n}L_{r}(n,k)(x)_{k}\big((x-k)+n+2r+k\big) \nonumber \\
	&=\ \sum_{k=0}^{n}L_{r}(n,k)(x)_{k+1}+\sum_{k=0}^{n}L(n,k)(x)_{k}(n+2r+k) \nonumber\\
	&=\ \sum_{k=1}^{n+1}L_{r}(n,k-1)(x)_{k}+\sum_{k=0}^{n}L_{r}(n,k)(n+2r+k)(x)_{k} \nonumber \\
	&=\ \sum_{k=0}^{n+1}\big(L_{r}(n,k-1)+(n+2r+k)L_{r}(n,k)\big)(x)_{k}. \nonumber
\end{align}
By Lemma 1, we also have 
\begin{equation}
	\langle x+2r\rangle_{n+1}=\sum_{k=0}^{n+1}L_{r}(n+1,k)
(x)_{k}. \label{19}
\end{equation}
By \eqref{18} and \eqref{19}, we obtain the following proposition. 
\begin{proposition}
	For $n,k,r\ge 0$ with $n\ge k$, we have 
	\begin{displaymath}
		L_{r}(n+1,k)=L_{r}(n,k-1)+(n+2r+k)L_{r}(n,k). 
	\end{displaymath}
\end{proposition}
From \eqref{14}, we derive the generating function of the $r$-extended Lah-Bell numbers as follows:
\begin{align}
\sum_{n=0}^{\infty}B_{n,r}^{L}\frac{t^{n}}{n!}\ &=\ \sum_{n=0}^{\infty}\bigg(\sum_{k=0}^{n}L_{r}(n,k)\bigg)\frac{t^{n}}{n!}\ =\ \sum_{k=0}^{\infty}\sum_{n=k}^{\infty}L_{r}(n,k)\frac{t^{n}}{n!} \label{21}\\
&=\ \sum_{k=0}^{\infty}\frac{1}{k!}\bigg(\frac{1}{1-t}-1\bigg)^{k}\bigg(\frac{1}{1-t}\bigg)^{2r}\ =\ e^{\big(\frac{1}{1-t}-1\big)}\bigg(\frac{1}{1-t}\bigg)^{2r}. \nonumber 	
\end{align}
Therefore, we obtain the following theorem. 
\begin{theorem}
	For $r\ge 0$, we have 
	\begin{displaymath}
		e^{\big(\frac{1}{1-t}-1\big)}\bigg(\frac{1}{1-t}\bigg)^{2r}= \sum_{n=0}^{\infty}B_{n,r}^{L}\frac{t^{n}}{n!}.
	\end{displaymath}
\end{theorem}
By Theorem 3, we get 
\begin{align}
	\sum_{n=0}^{\infty}B_{n,r}^{L}\frac{t^{n}}{n!}\ &=\ \frac{1}{e}e^{\frac{1}{1-t}}\cdot\bigg(\frac{1}{1-t}\bigg)^{2r}\ =\ \frac{1}{e}\sum_{k=0}^{\infty}\frac{1}{k!}\bigg(\frac{1}{1-t}\bigg)^{k+2r} \nonumber\\
	&=\ \sum_{n=0}^{\infty}\bigg(\frac{1}{e}\sum_{k=0}^{\infty}\frac{\langle k+2r\rangle_{n}}{k!}\bigg)\frac{t^{n}}{n!}. \label{22}
\end{align}
Therefore, by comparing the coefficients on both sides of \eqref{22}, we obtain the following theorem. 
\begin{theorem}[Dobinski-like formula]
	For $n\ge 0$, we have 
	\begin{displaymath}
		B_{n,r}^{L}=\frac{1}{e}\sum_{k=0}^{\infty}\frac{\langle k+2r\rangle_{n}}{k!}. 
	\end{displaymath}
\end{theorem}
In view of \eqref{2}, we define {\it{the $r$-extended Lah-Bell polynomials as}}
\begin{equation}
	B_{n,r}^{L}(x)=\sum_{k=0}^{n}L_{r}(n,k)x^{k},\quad(n\ge 0). \label{23}
\end{equation}
From \eqref{23}, we note that 
\begin{align}
\sum_{n=0}^{\infty}B_{n,r}^{L}(x)\frac{t^{n}}{n!}&=\ e^{x\big(\frac{1}{1-t}-1\big)}\bigg(\frac{1}{1-t}\bigg)^{2r} \label{24} \\
&=\ \frac{1}{e^{x}}\sum_{k=0}^{\infty}\bigg(\frac{1}{1-t}\bigg)^{k+2r}\frac{x^{k}}{k!}\nonumber \\
&=\ \frac{1}{e^{x}}\sum_{k=0}^{\infty}\frac{x^{k}}{k!}\sum_{n=0}^{\infty}\langle k+2r\rangle_{n}\frac{t^{n}}{n!}\nonumber \\
&=\ \sum_{n=0}^{\infty}\bigg\{\frac{1}{e^{x}}\sum_{k=0}^{\infty}\frac{\langle k+2r\rangle_{n}}{k!}x^{k} \bigg\}\frac{t^n}{n!}. \nonumber
\end{align}
Therefore, by \eqref{24}, we obtain the following theorem. 
\begin{theorem}
For $n\ge 0$, we have 
\begin{displaymath}
B_{n,r}^{L}(x)=\frac{1}{e^{x}}\sum_{k=0}^{\infty}\frac{\langle k+2r\rangle_{n}}{k!}x^{k}. 
\end{displaymath}
\end{theorem}
It is known that the $r$-extended Bell numbers are given by 
\begin{equation}
e^{e^{t}-1+rt}=\sum_{n=0}^{\infty}B_{n,r}\frac{t^{n}}{n!},\quad(\mathrm{see}\ [9]). \label{25}
\end{equation} \par
Replacing $t$ by $-\log(1-t)$ in \eqref{25}, we get 
\begin{align}
e^{\frac{1}{1-t}-1}\bigg(\frac{1}{1-t}\bigg)^{2r}\ &=\ \sum_{k=0}^{\infty}B_{k,2r}\frac{1}{k!}\big(-\log(1-t)\big)^{k} \nonumber \\
&=\ \sum_{k=0}^{\infty}B_{k,2r}(-1)^{k}\sum_{n=k}^{\infty}(-1)^{n}S_{1}(n,k)\frac{t^{n}}{n!} \label{26} \\
&=\ \sum_{n=0}^{\infty}\bigg(\sum_{k=0}^{n}(-1)^{n-k}S_{1}(n,k)B_{k,2r}\bigg)\frac{t^{n}}{n!}. \nonumber
\end{align}
From Theorem 3 and \eqref{26}, we have 
\begin{equation}
	B_{n,r}^{L}=\sum_{k=0}^{n}(-1)^{n-k}S_{1}(n,k)B_{k,2r},\quad(n,r\ge 0). \label{27}
\end{equation} \par
Replacing $t$ by $1-e^{-t}$ in Theorem 3, we have
\begin{align}
e^{e^{t}-1+2rt}\ &=\ \sum_{k=0}^{\infty}B_{k,r}^{L}\frac{1}{k!}\big(1-e^{-t}\big)^{k} \nonumber \\
&=\ \sum_{k=0}^{\infty}(-1)^{k}B_{k,r}^{L}\sum_{n=k}^{\infty}S_{2}(n,k)(-1)^{n}\frac{t^{n}}{n!}\label{28} \\
&=\ \sum_{n=0}^{\infty}\bigg(\sum_{k=0}^{n}(-1)^{n-k}S_{2}(n,k)B_{k,r}^{L}\bigg)\frac{t^{n}}{n!}. \nonumber
\end{align}
By \eqref{25} and \eqref{28}, we get 
\begin{equation}
	B_{n,2r}=\sum_{k=0}^{n}(-1)^{n-k}S_{2}(n,k)B_{k,r}^{L},\quad(n,k\ge 0). \label{29}
\end{equation}
Therefore, we obtain the following theorem. 
\begin{theorem}
	For $n\ge 0$, we have 
	\begin{displaymath}
		B_{n,r}^{L}=\sum_{k=0}^{n}(-1)^{n-k}S_{1}(n,k)B_{k,2r},
	\end{displaymath}
	and 
	\begin{displaymath}
		B_{n,2r}=\sum_{k=0}^{n}(-1)^{n-k}S_{2}(n,k)B_{k,r}^{L}. 
	\end{displaymath}
\end{theorem}
From \eqref{13}, we note that 
\begin{align}
\sum_{n=k}^{\infty}L_{r}(n,k)\frac{t^{n}}{n!}\ &=\ \frac{1}{k!}\bigg(\frac{t}{1-t}\bigg)^{k}\bigg(\frac{1}{1-t}\bigg)^{2r}\label{30}	 \\
&=\ \sum_{m=k}^{\infty}L(m,k)\frac{t^{m}}{m!}\sum_{l=0}^{\infty}\langle 2r\rangle_{l}\frac{t^{l}}{l!} \nonumber \\
&=\ \sum_{n=k}^{\infty}\bigg(\sum_{m=k}^{n}\binom{n}{m}L(m,k)\langle 2r\rangle_{n-m}\bigg)\frac{t^{n}}{n!}.\nonumber
\end{align}
Therefore, by \eqref{30}, we obtain the following lemma. 
\begin{lemma}
	For $n,k,r\ge 0$, with $n\ge k$, we have 
	\begin{displaymath}
		L_{r}(n,k)= \sum_{m=k}^{n}\binom{n}{m}L(m,k)\langle 2r\rangle_{n-m}.
	\end{displaymath}
\end{lemma}
Letting $x=2r$, and replacing $t$ by $-\log(1-t)$ in \eqref{8}, we have 
\begin{align}
\frac{1}{k!}\bigg(\frac{1}{1-t}\bigg)^{2r}\bigg(\frac{1}{1-t}-1\bigg)^{k}\ &=\ \sum_{m=k}^{\infty}S_{2}(m,k|2r)\frac{1}{m!}\big(-\log(1-t)\big)^{m}\label{31}\\
&=\ \sum_{m=k}^{\infty}(-1)^{m}S_{2}(m,k|2r)\sum_{n=m}^{\infty}(-1)^{n}S_{1}(n,m)\frac{t^{n}}{n!}\nonumber \\
&=\ \sum_{n=k}^{\infty}\bigg(\sum_{m=k}^{n}(-1)^{n-m}S_{2}(m,k|2r)S_{1}(n,m)\bigg)\frac{t^{n}}{n!}. \nonumber
\end{align}
Therefore, by \eqref{13} and \eqref{31}, we obtain the following theorem. 
\begin{theorem}
	For $n,k,r\ge 0$, with $n\ge k$, we have 
	\begin{displaymath}
		L_{r}(n,k)=\sum_{m=k}^{n}(-1)^{n-m}S_{2}(m,k|2r)S_{1}(n,m). 
	\end{displaymath}
\end{theorem}
Replacing $t$ by $1-e^{-t}$ in \eqref{13}, we get 
\begin{align}
\frac{1}{k!}\big(e^{t}-1\big)^{k}e^{2rt}\ &=\ \sum_{m=k}^{\infty}L_{r}(m,k)\frac{1}{m!}\big(1-e^{-t}\big)^{m} \nonumber \\
&=\ \sum_{m=k}^{\infty}(-1)^{m}L_{r}(m,k)\sum_{n=m}^{\infty}S_{2}(n,m)(-1)^{n}\frac{t^{n}}{n!} \label{32} \\
&=\ \sum_{n=k}^{\infty}\bigg(\sum_{m=k}^{n}(-1)^{n-m}S_{2}(n,m)L_{r}(m,k)\bigg)\frac{t^{n}}{n!}. \nonumber
\end{align}
Now, from \eqref{8} and \eqref{32}, we obtain the following corollary. 
\begin{corollary}
	For $n,k,r\ge 0$, with $n\ge k$, we have 
	\begin{displaymath}
		S_{2}(n,k|2r)=\sum_{m=k}^{n}(-1)^{m-m}S_{2}(n,m)L_{r}(m,k). 
	\end{displaymath}
\end{corollary}
For $m,n,k\in\mathbb{N}$, we observe from \eqref{13} that 
\begin{align}
\frac{1}{m!}&\bigg(\frac{1}{1-t}-1\bigg)^{m}\bigg(\frac{1}{1-t}\bigg)^{2r}\frac{1}{k!}\bigg(\frac{1}{1-t}-1\bigg)^{k}\bigg(\frac{1}{1-t}\bigg)^{2r} \label{34} \\
&=\ \frac{(m+k)!}{m!k!}\frac{\big(\frac{1}{1-t}-1\big)^{m+k}}{(m+k)!}\bigg(\frac{1}{1-t}\bigg)^{4r}\nonumber \\
&=\ \binom{m+k}{m}\sum_{n=m+k}^{\infty}L_{2r}(n,m+k)\frac{t^{n}}{n!}. \nonumber
\end{align}
On the other hand,
\begin{align}
&\frac{1}{m!}\bigg(\frac{1}{1-t}-1\bigg)^{m}\bigg(\frac{1}{1-t}\bigg)^{2r}\frac{1}{k!}\bigg(\frac{1}{1-t}-1\bigg)^{k}\bigg(\frac{1}{1-t}\bigg)^{2r} \label{35}\\
&\ =\ \sum_{l=m}^{\infty}L_{r}(l,m)\frac{t^{l}}{l!}\sum_{j=k}^{\infty}L_{r}(j,k)\frac{t^{j}}{j!} \nonumber \\
&\ =\ \sum_{n=m+k}^{\infty}\bigg(\sum_{l=m}^{n-k}\binom{n}{l}L_{r}(l,m)L_{r}(n-l,k)\bigg)\frac{t^{n}}{n!}\nonumber.
\end{align}
Therefore, by \eqref{34} and \eqref{35} we obtain the following theorem. 
\begin{theorem}
	For $m,n,k\in\mathbb{N}$, with $n\ge m+k$ and $r\ge 0$, we have 
	\begin{displaymath}
		\binom{m+k}{m}L_{2r}(n,m+k)=\sum_{l=m}^{n-k}\binom{n}{l}L_{r}(l,m)L_{r}(n-l,k). 
	\end{displaymath}
\end{theorem}
In particular, we have 
\begin{displaymath}
\binom{m+k}{m}L(n,m+k)=\sum_{l=m}^{n-k}L(l,m)L(n-l,k).
\end{displaymath}

\section{Further Remarks} 
Let $X$ be the Poission random variable with parameter $\alpha(>0)$. 
From \eqref{11}, we note that 
\begin{align}
E\bigg[\bigg(\frac{1}{1-t}\bigg)^{X}\bigg]\ &=\ \sum_{i=0}^{\infty}\bigg(\frac{1}{1-t}\bigg)^{i}e^{-\alpha}\frac{\alpha^{i}}{i!}\label{36} \\
&=e^{\alpha\big(\frac{1}{1-t}-1\big)}\ =\ \sum_{n=0}^{\infty}B_{n}^{L}(\alpha)\frac{t^{n}}{n!}.\nonumber 
\end{align}
On the other hand, 
\begin{equation}
E\bigg[\bigg(\frac{1}{1-t}\bigg)^{X}\bigg]=\sum_{n=0}^{\infty}E[\langle X\rangle_{n}]\frac{t^{n}}{n!}. \label{37}
\end{equation}
Thus, we see that 
\begin{displaymath}
E[\langle X\rangle_{n}]=B_{n}^{L}(\alpha),\quad(n\ge 0). 
\end{displaymath} \par
Now, we observe from \eqref{24} that 
\begin{align}
E\bigg[\bigg(\frac{1}{1-t}\bigg)^{X+2r}\bigg]\ &=\ \sum_{k=0}^{\infty}\bigg(\frac{1}{1-t}\bigg)^{k}e^{-\alpha}\frac{\alpha^{k}}{k!}\bigg(\frac{1}{1-t}\bigg)^{2r} \label{38} \\
&=\ e^{\alpha\big(\frac{1}{1-t}-1\big)}\bigg(\frac{1}{1-t}\bigg)^{2r} \nonumber \\
&=\ \sum_{n=0}^{\infty}B_{n,r}^{L}(\alpha)\frac{t^{n}}{n!}.\nonumber
\end{align}
On the other hand, 
\begin{equation}
	E\bigg[\bigg(\frac{1}{1-t}\bigg)^{X+2r}\bigg]=\sum_{n=0}^{\infty}E[\langle X+2r\rangle_{n}]\frac{t^{n}}{n!}. \label{39}
\end{equation}
Therefore, by \eqref{38} and \eqref{39}, we obtain the following theorem. 
\begin{theorem}
Let $X$ be the Poisson random variable with parameter $\alpha(>0)$. Then we have 
\begin{displaymath}
E[\langle X+2r\rangle_{n}]=B_{n,r}^{L}(\alpha),\quad(n\ge 0).
\end{displaymath}
\end{theorem}
From Lemma 1, we note that 
\begin{align*}
	\langle X+2r\rangle_{n}\ &=\ \sum_{k=0}^{n}L_{r}(n,k)(X)_{k} \ =\ \sum_{k=0}^{n}\sum_{l=0}^{k}L_{r}(n,k)S_{1}(k,l)X^{l} \\
	&=\ \sum_{l=0}^{n}\sum_{k=l}^{n}L_{r}(n,k)S_{1}(k,l)X^{l}.
\end{align*}
Thus, we note that 
\begin{equation}
B_{n,r}^{L}(\alpha)=E[\langle X+2r\rangle_{n}]=\sum_{l=0}^{n}\sum_{k=l}^{n}L_{r}(n,k)S_{1}(k,l)E[X^{l}]. \label{40}
\end{equation}
By \eqref{9}, we get 
\begin{equation}
	E[X^{l}]=\sum_{k=0}^{\infty}k^{l}p(k)=e^{-\alpha}\sum_{k=0}^{\infty}\frac{k^{l}}{k!}\alpha^{k}=B_{l}(\alpha),\quad(l\ge 0). \label{41}
\end{equation}
From \eqref{40} and \eqref{41}, we have 
\begin{displaymath}
	B_{n,r}^{L}(\alpha)=\sum_{l=0}^{n}\sum_{k=l}^{n}L_{r}(n,k)S_{1}(k,l)B_{l}(\alpha),\quad(n\ge 0). 
\end{displaymath}
Let $X$ be the Poisson random variable with parameter $\alpha>0$. Then we have 
\begin{displaymath}
	B_{n,r}^{L}(\alpha)=E[\langle X+2r\rangle_{n}]=\sum_{l=0}^{n}\sum_{k=l}^{n}L_{r}(n,k)S_{1}(k,l)B_{l}(\alpha). 
\end{displaymath}

\section{Conclusion}
In this paper, we obtained a recurrence relation for $r$-Lah numbers and the generating function of the $r$-extended Lah-Bell numbers. Also, we derived Dobinski-like formulas for the $r$-extended Lah-Bell numbers and polynomials. We expressed the $r$-extended Lah-Bell numbers in terms of the 
$r$-extended Bell numbers and vice versa. Also, we deduced an expression of $r$-Lah numbers in terms of the Stirling polynomials of the second kind and vice versa. In addition, we were able to get a relation among $r$-Lah numbers. Lastly, we showed that the rising factorial moments of $X+2r$ is equal to the $r$-extended Lah-Bell numbers evaluated at $\alpha$, where $X$ is the Poisson random variable with parameter $\alpha (>0)$, and $r$ is a nonnegative integer. This follows from the fact that the probability-generating function of the Poisson random variable shifted by an even nonnegative integer $2r$ is equal to the generating function of $r$-exteneded Lah-Bell polynomials. \par 
It is Carlitz who initiated a study of degenerate Bernoulli and Euler numbers and polynomials. Following his pioneering work,  in recent years we have explored some special numbers and polynomials and their degenerate versions, and discovered their arithmetical and combinatorial properties and some of their applications. It is one of our furure projects to continue to work on these by exploiting various means such as generating functions, combinatorial methods, $p$-adic analysis, umbral calculus, differential equations and probability theory. \par

\end{document}